\documentclass[11pt]{amsart}
\usepackage{url}
\usepackage{bm}
\usepackage[left=1in,top=1in,right=1in,bottom=1in,head=.2in]{geometry}

\usepackage{enumitem}
\usepackage[textsize=scriptsize]{todonotes}

\usepackage{fancyhdr}
\newtheorem*{theorem*}{Theorem}

\newtheorem*{lemma*}{Lemma}

\newtheorem*{corollary*}{Corollary}
\theoremstyle{definition}

\newtheorem*{examples*}{Examples}
\newtheorem*{remark*}{Remark}
\newtheorem*{remarks*}{Remarks}
\newtheorem*{addenda*}{Addenda}

\newcommand{\figref}[1]{Figure~\ref{F:#1}}

\newcommand{\fig}[3]{\begin{figure}[h!] \includegraphics[height=#1pt]{#2}#3\end{figure}}

\newcommand{\bc}{\mathbb C}

\newcommand{\sss}{S^2\!\times\!S^2}
\newcommand{\sts}{S^2\widetilde\times S^2}
\newcommand{\cs}{\,\#\,}

\newcommand{\cptwo}{\bc P^2}
\newcommand{\cptwobar}{\smash{\overline{\bc P}^2}}

\newcommand{\del}{\partial}
\newcommand{\s}{\Sigma}

\newcommand{\items}{\begin{itemize}[leftmargin=25pt,rightmargin=5pt]
  \setlength\itemsep{2pt}}
\newcommand{\stopitems}{\end{itemize}}

\begin{document}

\title{Isotopy of surfaces in $4$-manifolds \\ after a single stabilization}
\author[Dave Auckly]{Dave Auckly$^{1}$}
\address{Department of Mathematics\newline\indent Kansas State University\newline\indent  Manhattan,
Kansas 66506}
\email{dav@math.ksu.edu}
\author[Hee Jung Kim]{Hee Jung Kim$^{1,2}$}
\address{Department of Mathematical Sciences\newline\indent
Seoul National University\newline\indent
1 Gwanak-ro\newline\indent
Seoul, 08826, South Korea}
\email{heejungorama@gmail.com}
\author[Paul Melvin]{Paul Melvin$^{1}$}
\address{Department of Mathematics\newline\indent  Bryn Mawr College\newline\indent  Bryn Mawr, PA 19010}
\email{pmelvin@brynmawr.edu}
\author[Daniel Ruberman]{Daniel Ruberman$^{1,3}$}
\address{Department of Mathematics, MS 050\newline\indent Brandeis
University \newline\indent Waltham, MA 02454}
\email{ruberman@brandeis.edu}
\author[Hannah Schwartz]{Hannah Schwartz$^{1}$}
\address{Department of Mathematics\newline\indent  Bryn Mawr College\newline\indent  Bryn Mawr, PA 19010}
\email{hrschwartz@brynmawr.edu}
\thanks{$^1$All of the authors were supported by an AIM SQuaRE grant.
$^2$Supported by  NRF grants 2015R1D1A1A01059318 and BK21 PLUS SNU Mathematical Sciences Division. 
$^3$Partially supported by NSF Grant 1506328.
}

\maketitle

\vskip-25pt
\vskip-25pt

\begin{abstract}
Any two homologous surfaces of the same genus embedded in a smooth 4-manifold $X$ with simply-connected complements are shown to be smoothly isotopic in $X\cs\sss$ if the surfaces are ordinary, and in $X\cs\sts$ if they are characteristic. 
\end{abstract}

\parskip 3pt

\section{Introduction}

By the $4$-dimensional light bulb theorem of Gabai \cite[Theorem 1.2]{gabai}, any two homologous $2\text{-spheres}$ embedded with a common geometric dual in a smooth simply-connected $4\text{-manifold}$ $X$ are isotopic; here a {\it geometric dual} for a surface $F\subset X$ is an embedded $2$-sphere $\s$ of {\it square} zero (meaning self-intersection zero) intersecting $F$ transversely in a single point.  The same result holds for homologous closed surfaces $F_0$ and $F_1$ of the same genus embedded in $X$ under a mild fundamental group condition (that the $F_i$ should be ``$\s$--inessential"; see \cite[Theorem 9.7]{gabai}).  

It has been known for some time that this result fails without a common geometric dual for the surfaces $F_i$, or even without the weaker condition that each $F_i$ should have an {\sl immersed} geometric dual, i.e.\ a simply connected complement.  Examples arise easily from the existence of exotic smooth structures on closed simply-connected $4$-manifolds (Donaldson \cite{donaldson}) and the fact that such manifolds become diffeomorphic after sufficiently many stabilizations (Wall \cite{wall:4manifolds}); here a {\it stabilization} means a connected sum with $\sss$.  Work of Quinn \cite{quinn} and Perron \cite{perron1,perron2} shows that the surfaces $F_i$ always become isotopic after sufficiently many {\it external} stabilizations, where the connected sums are taken away from $F_0\cup F_1$.  This raises the question of how many stabilizations are needed.  In particular, {\it is one enough}?  If $n$ is the minimal number of stabilizations needed, the surfaces are said to be {\it strictly $n$-stably isotopic}; the first explicit examples of 
families of strictly $1$-stably isotopic surfaces were given by the authors in \cite{akmr} (see also Akbulut \cite{akbulut}).  

Analogous questions have been asked about many exotic phenomena in 4-dimensional topology which are known to dissipate after sufficiently many stabilizations.  To the authors' knowledge, there are no instances known where it can be shown that one is {\sl not} enough.  

In this note it is shown using Gabai's results that, indeed, one is always enough if the surfaces have simply-connected complements and are {\it ordinary}.  Recall that a surface is ordinary if it is not {\it characteristic}, meaning dual to the second Stiefel-Whitney class $w_2(X)$; geometrically, a surface is characteristic if it intersects even classes evenly and odd classes oddly.  Furthermore, this one-is-enough result still holds for {\it characteristic} surfaces under `twisted" stabilization, meaning connected sum with the twisted bundle $\sts$.         

\section{One is enough}

\begin{theorem*}\label{T:oneisenough}
If $X$ is a smooth simply-connected $4$-manifold and $\alpha\in H_2(X)$ is an ordinary class, 
then any two closed oriented surfaces $F_0$ and $F_1$ in $X$ of the same genus representing $\alpha$, both with simply-connected complement, are smoothly isotopic in $X\cs \sss$ $($summing away from $F_0\cup F_1)$.  When $\alpha$ is characteristic, the same result holds if one stabilizes by summing with $\sts$.
\end{theorem*}

\proof  It can be assumed that $F_0$ and $F_1$ intersect transversely.  When $\alpha$ is ordinary, the strategy is to find an embedded sphere in $X\cs(\sss)$ of square zero geometrically dual to both $F_0$ and $F_1$.  The result will then follow from Gabai's theorem.  When $\alpha$ is characteristic, no such sphere exists since the image of $\alpha$ under the natural map $H_2(X) \to H_2(X\cs\sss)$ is still characteristic.  But it will be seen to exist in $X\cs\sts$, and the proof will follow as before.  Here is how one carries out this strategy, assuming first that $\alpha$ is ordinary.   
 
Start with an {\sl immersed} sphere $\s\subset X$ meeting $F_0$ transversally in one point, i.e.\ an immersed geometric dual for $F_0$, that is transverse to $F_1$.  Such a sphere exists since $\pi_1(X-F_0)=1$. 
If $X$ is even, then $\s$ has even square.  If $X$ is odd, then $\s$ may have odd square, but it can be modified to have even square by connected summing with an immersed sphere $S$ in $X-F_0$ of odd square.  To see that such an $S$ exists, note that since $F_0$ is ordinary, there exists an immersed surface $E \subset X$ whose square is of the opposite parity from the algebraic intersection number $n=E\cdot F_0$.  Now by the immersed Norman trick~\cite{norman:trick}, $E$ can be tubed to parallel copies of $\s$ along arcs in $F_0$ to remove its intersections with $F_0$, giving a surface $E+n\s$ of {\sl odd} square in $X-F_0$, since $E+n\s$ is homologically the sum of an even and odd class.  But since $X-F_0$ is simply-connected, this last surface is homologous to an immersed sphere $S$ in $X-F_0$.  Thus it can always be arranged for the self-intersection of $\s$ to be even.

Orient $\s$ so that $\s\cdot F_0 = 1$.  Then $\s\cdot F_1 = 1$ as well, since $F_1$ is homologous to $F_0$.  If the geometric intersection number $|\s \cap F_1|$ is greater than $1$, then finger and Whitney moves 
can be used to reduce this number, thus inductively moving $\s$ to meet $F_1$ in only one point.  This is accomplished as follows:  

For any oppositely oriented pair $p,\,q$ of intersection points in $\s\cap F_1$, consider a Whitney circle on $\s\cup F_1$ disjoint from $F_0$, made up of two arcs $\gamma_0\subset \s$ and $\gamma_1\subset F_1$ meeting at their endpoints $p$ and $q$.  
Since $X-F_1$ is simply-connected, this circle bounds an immersed disk $D$ with interior disjoint from $F_1$, but not necessarily disjoint from $F_0$ and $\Sigma$.  See \figref{config} for a schematic of the intersections and self-intersections between $F_0$, $F_1$, $\s$ and $D$.    

\fig{150}{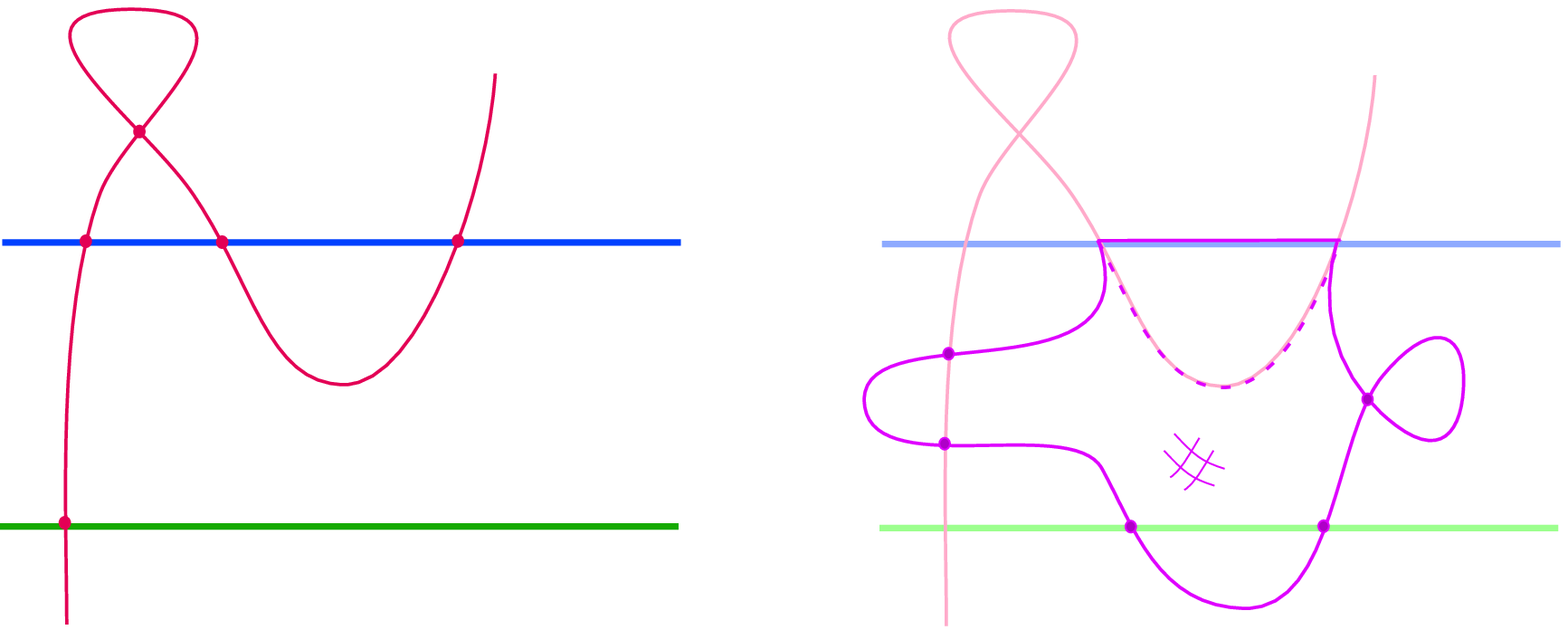}{
\put(-250,12){$F_0$}
\put(-250,80){$F_1$}
\put(-375,55){$\s$}
\put(-20,55){$D$}
\put(-113,100){\small$p$}
\put(-58,100){\small$q$}
\put(-88,102){$\gamma_1$}
\put(-88,63){$\gamma_0$}
\caption{The surfaces $F_0$ and $F_1$, the immersed dual sphere $\s$, and a Whitney disk $D$}
\label{F:config}}

\noindent Now perform finger moves of $F_0$ across $F_1$ to remove the intersections of $F_0$ with $D$, guided by disjoint arcs in $D-\s$ from the points of $D\cap F_0$ to $\gamma_1$.  To perform a Whitney move of $\s$ across $D$ one must first fix the framing, that is, ensure that the restriction to $\del D$ of the framing of the normal bundle of $D$ matches the framing induced by $\s\cup F_1$.  This is achieved by ``boundary twisting" $D$ around $\s$ along $\gamma_0$ (see Freedman and Quinn \cite[\S1.3--1.4]{freedman-quinn} for more details). This introduces additional intersection points between $\s$ and $D$, but keeps the interior of $D$ and $F_0 \cup F_1$ disjoint.  
Once the framing has been fixed, one may push $\s$ over the immersed Whitney disk $D$, thus removing its intersection points $p$ and $q$ with $F_1$. Repeating this process, one obtains an immersed sphere $\s\subset X$ with even self-intersection that is geometrically dual to both $F_0$ and $F_1$.

To arrange for $\s$ to be embedded, the ambient manifold $X$ must be stabilized.  In particular, take the connected sum with $\sss$ at a point in the complement of $F_0\cup F_1\cup\s$. Choose coordinate 2-spheres $S = S^2\times \{\text{pt}\}$ and $T =  \{\text{pt}\}\times S^2$ in $\sss$.  Now replace $\s$ by its connected sum with $S$ along a tube disjoint from $F_0\cup F_1$, so that $\s$ now has an embedded geometric dual sphere $T$ in $X\cs\sss$.  Then eliminate the double points in $\s$ by the Norman trick, tubing to parallel copies of $T$.  The result is an {\sl embedded} sphere, still denoted $\s$, that intersects both $F_0$ and $F_1$ geometrically in exactly one point. This sphere still has even square and an embedded geometric dual sphere $T$ of square zero. Tubing with additional copies of $T$, the dual sphere $\s$ can be made to have self-intersection zero. 

This entire argument can be repeated when $\alpha$ is characteristic, except that the self-intersection of the immersed dual $\s$ for $F_0$ and $F_1$ will  now necessarily be odd.  In this case, stabilize by summing with $\sts$, viewed as a Hirzebruch surface with a section $S$ of odd square, and fiber $T$. Now tube $\s$ to $S$ to obtain an immersed sphere of even square and a geometric dual $T$, and proceed as in the ordinary case by tubing with $T$ to make $\s$ embedded.

The proof is finished by applying Gabai's results \cite{gabai}.  The technical assumption needed for surfaces of higher genus (that $\pi_1(F_i-\s) \to\pi_1(Y - \s)$ for $F_i\subset Y$ should be trivial; see Theorem 9.7 in \cite{gabai})
holds since $\s$ has a geometric dual $T$, so the complement of $\s$ in $Y= X\cs\sss$ (or $Y=X\cs\sts$) is simply-connected.  \qed

\medskip

As a consequence of this theorem, one can find infinite families of 1-stably equivalent 2-spheres in many once-stabilized 4-manifolds. 


\begin{corollary*}
Let $X_1,X_2\dots$ be any $($possibly infinite$)$ list of pairwise non-diffeomorphic, smooth, closed, simply-connected $4$-manifolds, all homeomorphic to one such $X$, such that $X_i\cs\sss$ is diffeomorphic to $X\cs\sss$ for each $i$.  Also assume that the $X_i$ remain distinct after connected summing with any number of copies of $\cptwobar$.  If either 
\items
\item[\bf 1)] $X$ is even and $n$ is any even nonnegative integer, or
\item[\bf 2)] $X$ is odd and indefinite, and $n$ is any nonnegative integer, 
\stopitems
then there is a corresponding family of strictly $1$-stably isotopic 2-spheres $S_1,S_2,\dots$ of square $n$ smoothly embedded in $X\cs\sss$.  
\end{corollary*}

\proof
First observe that if $X$ is even, then it is indefinite by Donaldson's Theorem A \cite{donaldson:A}, and so in either case, every automorphism of the quadratic form of $X\cs\sss$ is induced by a diffeomorphism of $X\cs\sss$ \cite{wall:diffeomorphisms}.  

First assume that $n$ is even.  Then $S^2\times S^2$ is diffeomorphic to the $S^2$-bundle over $S^2$ of Euler class $n$.  Let $S$ denote the zero section  of this bundle, which is an embedded 2-sphere of square $n$, and $T$ denote the fiber.  By abuse of notation, $S$ and $T$ will also denote the corresponding 2-spheres in the $\sss$ factor in $X\cs\sss$, and also in $X_i\cs\sss$ for each $i$.  By Wall's result, there exist diffeomorphisms $h_i\colon X_i\cs\sss \to X\cs\sss$ such that the 2-spheres $S_i=h_i(S)$ are all homologous to $S$ in $X\cs\sss$.  However, these spheres are not smoothly isotopic in $X\cs\sss$,  since blowing up $n$ points on $S_i$ and then surgering the resulting sphere yields $X_i\cs n\cptwobar$ (seen for example by a Kirby calculus exercise), and these manifolds are distinct by hypothesis. 

If $n$ is odd and $X$ is odd, use the facts that $X\cs\sss$ is diffeomorphic to $X\cs\sts$ and that $\sts$ is diffeomorphic to the $S^2$-bundle over $S^2$ of Euler class $n$.  As above, one then constructs a family of homologous but non-isotopic spheres $S_i\subset X\cs\sss$ of square $n$.


Now in either case, the spheres $S_i$ have geometric dual spheres $T_i=h_i(T)$ of self-intersection zero, so the $S_i$ are ordinary with simply-connected complements.  It follows from the theorem that the $S_i$ become isotopic in $X\cs\sss\cs\sss$.  Since the $S_i$ are not isotopic in $X\cs\sss$, they are strictly 1-stably isotopic.      \qed

\begin{examples*} {\bf 1)} Let $K_i$ be a sequence of knots with distinct Alexander polynomials, and $X_i$ be the 4-manifolds obtained from the elliptic surface $E(2)$ by $K_i$-knot surgery along a regular fiber.  Then the $X_i$ are pairwise non-diffeomorphic \cite{fintushel-stern, sunukjian} and satisfy the stability hypothesis of the corollary \cite{akbulut:stable,auckly}.  Thus the spin manifold $E(2)\cs\sss$ contains strictly 1-stably isotopic families of spheres of any even, nonnegative self-intersection.  

\noindent{\bf 2)} Let $X_i$ be the family of Dolgacev surfaces obtained from the rational elliptic surface $E(1)$ by a pair of logarithmic transforms of orders $2$ and $2i+1$.  Then the $X_i$ are pairwise non-diffeomorphic $4$-manifolds \cite{donaldson,friedman-morgan,okonek-ven:dolgachev} that satisfy the stability hypothesis of the corollary \cite{mandelbaum}.   Thus $2\cptwo\cs10\cptwobar$ contains strictly 1-stably isotopic families of spheres of any nonnegative self-intersection. 
\end{examples*}


\end{document}